\documentclass[preprint,1p,10pt]{elsarticle}

\usepackage{graphicx}
\usepackage{amssymb}

\usepackage{amsmath}
\usepackage{cleveref}

\usepackage{caption}

\usepackage{color}
\usepackage[table]{xcolor}

\biboptions{sort&compress,comma,round}

\journal{CMAME}

\begin{document}

\begin{frontmatter}


\title{Spectral approximation properties of isogeometric analysis with variable continuity}



\author[curtin,cic]{Vladimir Puzyrev\corref{cor1}}
\ead{vladimir.puzyrev@gmail.com}

\author[curtin,cic]{Quanling Deng}
\ead{qdeng12@gmail.com}

\author[curtin,cic,csiro]{Victor Calo}
\ead{vmcalo@gmail.com}

\cortext[cor1]{Corresponding author}

\address[curtin]{Department of Applied Geology,
Western Australian School of Mines, Curtin University,
Kent Street, Bentley, Perth, WA 6102, Australia}
\address[cic]{Curtin Institute for Computation,
Curtin University, Kent Street, Bentley, Perth, WA 6102, Australia}
\address[csiro]{Mineral Resources,
Commonwealth Scientific and Industrial Research Organisation (CSIRO),
Kensington, Perth, WA 6152, Australia}

\begin{abstract}
We study the spectral approximation properties of isogeometric analysis with local continuity reduction of the basis. Such continuity reduction results in a reduction in the interconnection between the degrees of freedom of the mesh, which allows for large savings in computational requirements during the solution of the resulting linear system. The continuity reduction results in extra degrees of freedom that modify the approximation properties of the method. The convergence rate of such refined isogeometric analysis is equivalent to that of the maximum continuity basis. We show how the breaks in continuity and inhomogeneity of the basis lead to artefacts in the frequency spectra, such as stopping bands and outliers, and present a unified description of these effects in finite element method, isogeometric analysis, and refined isogeometric analysis. Accuracy of the refined isogeometric analysis approximations can be improved by using non-standard quadrature rules. In particular, optimal quadrature rules lead to large reductions in the eigenvalue errors and yield two extra orders of convergence similar to classical isogeometric analysis.
\end{abstract}


\end{frontmatter}


\section{Introduction}
\label{S:1}

Isogeometric analysis (IGA) is a numerical technique for approximating the solutions of partial differential equations, which was introduced in 2005 \citep{hughes2005isogeometric} and received significant attention since then \citep{bazilevs2010isogeometric, cottrell2006isogeometric, cottrell2007studies, cottrell2009isogeometric, hughes2008duality, hughes2010efficient, hughes2014finite, gomez2008isogeometric, calo2008simulation, auricchio2013locking}. Spectrum analysis of isogeometric discretizations shows that this method is more accurate compared to the classical finite element analysis (FEA) for a fixed number of degrees of freedom \cite{cottrell2007studies, hughes2008duality, hughes2014finite}. Isogeometric analysis uses as basis functions those employed in computer aided design (CAD) systems that can represent exactly many complex geometries relevant in engineering applications. The isogeometric framework allows higher continuity across element interfaces (up to $p-1$ continuous derivatives across element boundaries, where $p$ is the polynomial order) and local control of the continuity of the basis, which is one of its most powerful and fundamental features \citep{cottrell2007studies}.

The choice of the numerical method heavily influences the performance of sparse linear solvers. For traditional FEA, the interconnection between the subdomains is weak due to the minimal inter-element continuity and they are connected by narrow separators. In highly continuous IGA, the increased support of basis functions strengthens the interconnection between the subdomains which results in wider separators. This leads to degradation of performance of both direct and iterative solvers for highly continuous IGA discretizations, increasing both time and memory requirements \citep{collier2012cost, collier2013cost, collier2014computational}. Recently, \citet{garcia2017value} proposed a modification of the continuity of the discrete space to exploit the increased continuity locally, while reducing the overall solution cost of the method. They called this optimized continuity discretization the refined isogeometric analysis (rIGA). This method reduces the continuity of the basis to achieve the high accuracy of classical IGA at a lower computational cost of the solution of linear systems than both IGA and FEA methods. In the simplest strategy, rIGA starts with the maximum $C^{p-1}$ continuity and reduces the continuity at certain hyperplanes in the mesh that act as separators during the elimination of the degrees of freedom employed by the direct solvers. This reduction in the interconnection between the degrees of freedom at critical zones allows for large savings in the computational cost which ultimately results in a reduction in the wall clock time needed to solve the resulting sparse systems of linear equations. rIGA performs localized reductions of continuity to optimize the computational complexity of the resulting isogeometric discretization. Each lower continuity separator results in a small increase in the total number of degrees of freedom in the system, while the total computational time of the solution decreases by a factor proportional to $p^2$ with respect to the isogeometric elements of the maximum continuity. This methodology reduces the computational complexity of the solution and the memory requirements for the solution. The total execution time can be reduced about fifty times in 3D. \citet{garcia2017optimally} reduce the continuity further between separators to achieve further complexity and memory usage savings.

We seek to quantify the approximation errors of a numerical method using dispersion and spectral analyses. The dispersive properties of classical FEA are well established in the literature, cf., \citep{thompson1994complex, ainsworth2004discrete} and references therein. The dispersion analysis of isogeometric elements revealed them to be superior to classical finite elements in approximating various types of partial differential equations \citep{bazilevs2007variational, hughes2008duality, hughes2014finite, kolman2014complex, dede2015isogeometric}. \citet{hughes2008duality} unified the dispersion and spectrum analysis showing their equivalence in the regime where the wavenumber is real. However, only uniform isogeometric elements of the highest $C^{p-1}$ continuity have been considered in most of these studies. The work of \citet{kolman2014complex} is perhaps the only study that considers spectrum analysis of B-spline multi-segment discretizations with $C^0$ continuity between these segments. This problem requires a comprehensive analysis since breaks in the continuity of the basis affect the approximation properties of the numerical method.

In this paper, we perform a unified dispersion and spectrum analysis of refined isogeometric analysis and study how the variable continuity of the basis affects the approximation quality. We address linear boundary- and initial-value problems by expressing them in terms of the eigenvalue and eigenfunction errors of the corresponding eigenproblem. We conduct a global error analysis, that is, we characterize the errors in the eigenvalues and the eigenfunctions for all the modes \citep{hughes2014finite}. Using Strang's Pythagorean eigenvalue error theorem \citep{strang1973analysis} and its generalization \citep{puzyrev2017dispersion}, we describe the total error budget which consists of the eigenvalue and eigenfunction errors. 

The outline of this work is as follows. First, we briefly describe the isogeometric framework and model problem under study in Sections 2 and 3. In Section 4, we analyze the approximation errors of rIGA and discuss its connections to standard isogeometric and finite elements. We analyze the stopping bands in Section 5, while we briefly describe the outliers in Section 6. Section 7 presents several possible ways to improve the accuracy of the rIGA approximations. Finally, in Section 8 we summarize our findings and discuss future research directions.

\section{B-spline basis functions and knot insertion: the isogeometric framework}
\label{S:2}

B-spline and Non-Uniform Rational B-spline (NURBS) basis functions \citep{piegl2012nurbs} are the most common geometrical representations used in the IGA framework, though other splines, such as T-splines, hierarchical splines and locally refined splines (LR-splines) are becoming common \citep{bazilevs2010isogeometric, schillinger2012isogeometric, dokken2013polynomial}. In this section, we briefly describe the B-spline basis functions and related concepts. For a more detailed
review, we refer the reader to \citep{piegl2012nurbs, hughes2005isogeometric, cottrell2009isogeometric}.

B-splines are piecewise polynomial curves composed of linear combinations of B-spline basis functions. For the piecewise constants ($p = 0$), they are defined as \citep{piegl2012nurbs}

\begin{equation} \label{eq:basis1}
N_{i,0}(\xi) = \left \{
\begin{array}{ll}
  1 \ \ \text{if} \ \xi_i \le \xi < \xi_{i+1}\\
  0 \ \ \text{otherwise}
  \end{array}
\right.
\ \ \ \ \ i=0, \dotsc, n,
\end{equation}
where $n$ is the number of basis functions which comprise the B-spline. For $p \ge 1$, the B-spline basis functions are defined recursively as

\begin{equation} \label{eq:basis2}
N_{i,p}(\xi) = \frac{\xi-\xi_i}{\xi_{i+p}-\xi_i} N_{i,p-1}(\xi) + \frac{\xi_{i+p+1}-\xi}{\xi_{i+p+1}-\xi_{i+1}} N_{i+1,p-1}(\xi).
\end{equation}

The set of non-decreasing real numbers $\xi_1, \dotsc, \xi_{n+p+1}$ represents the coordinates in the parametric space of the curve and is called the knot vector. The knots partition the parameter space into elements. Knot values may be repeated and these multiplicities of knot values have important implications for the continuity properties of the basis \citep{cottrell2007studies}. A knot vector is uniform if the knots are uniformly spaced and non-uniform otherwise. The open knots have the first and last knots repeated $p + 1$ times. B-Spline basis functions are non-negative and, for an open knot vector, constitute a partition of unity $\sum_{i=1}^N N_{i,p}(\xi) = 1 \ \ \forall \xi $. They are $C^{p-1}$-continuous inside the domain when internal knots are not repeated. For a knot of multiplicity $k \leq p$, the basis is $C^{p-k}$-continuous at that knot. For example, a $C^0$ separator can be added with the help of a knot of multiplicity $p$, thus making the basis to be $C^0$ at that knot.

By taking a linear combination of basis functions and control points, a piecewise polynomial B-Spline curve can be constructed. B-Spline surfaces and solids are constructed by means of tensor products \citep{cottrell2006isogeometric}. More details on B-splines, NURBS, and common geometric algorithms can be found in the book by \citet{piegl2012nurbs}.

The framework of B-splines allows for five tools for their management: knot insertion, knot refinement, knot removal, degree elevation, and degree reduction. Isogeometric analysis offers three different mechanisms for element refinement. In addition to the \textit{h}- and \textit{p}-refinement of the finite element framework, IGA allows for the \textit{k}-refinement \citep{cottrell2007studies}. Below we briefly review these three refinement techniques.

The mechanism for implementing \textit{h}-refinement in IGA is knot insertion. The insertion of \textit{new} knot values is similar to the classical \textit{h}-refinement in FEA as it reduces the support of the basis functions. Repeating \textit{existing} knot values decreases the continuity of the bases and their support, but does not reduce the element size. This mechanism, which is at the core of rIGA, has no analogue in standard finite elements.

\textit{p}-refinement is implemented in IGA using degree elevation. This mechanism involves increasing the polynomial order of the basis functions used to represent the geometry and the solution space. The multiplicity of each existing knot value is increased by one, but no new knot values are added. Pure \textit{p}-refinement increases the polynomial order while the basis remains $C^0$, thus leading to standard finite elements.

\textit{k}-refinement of IGA has no analogue in FEA. Pure \textit{k}-refinement keeps the mesh size fixed but increases the continuity along with the polynomial order. The solutions approximated with isogeometric elements of order $p$ have global continuity of order up to $p-1$, i.e., $p-1$ continuous derivatives across element boundaries. Highly smooth $C^{p-1}$ continuous basis functions produce better approximations of the derivatives of the solution when compared to $C^0$ finite elements for problems with smooth coefficients \citep{hughes2008duality}. 

\begin{figure}[!ht]
\centering\includegraphics[width=1.0\linewidth]{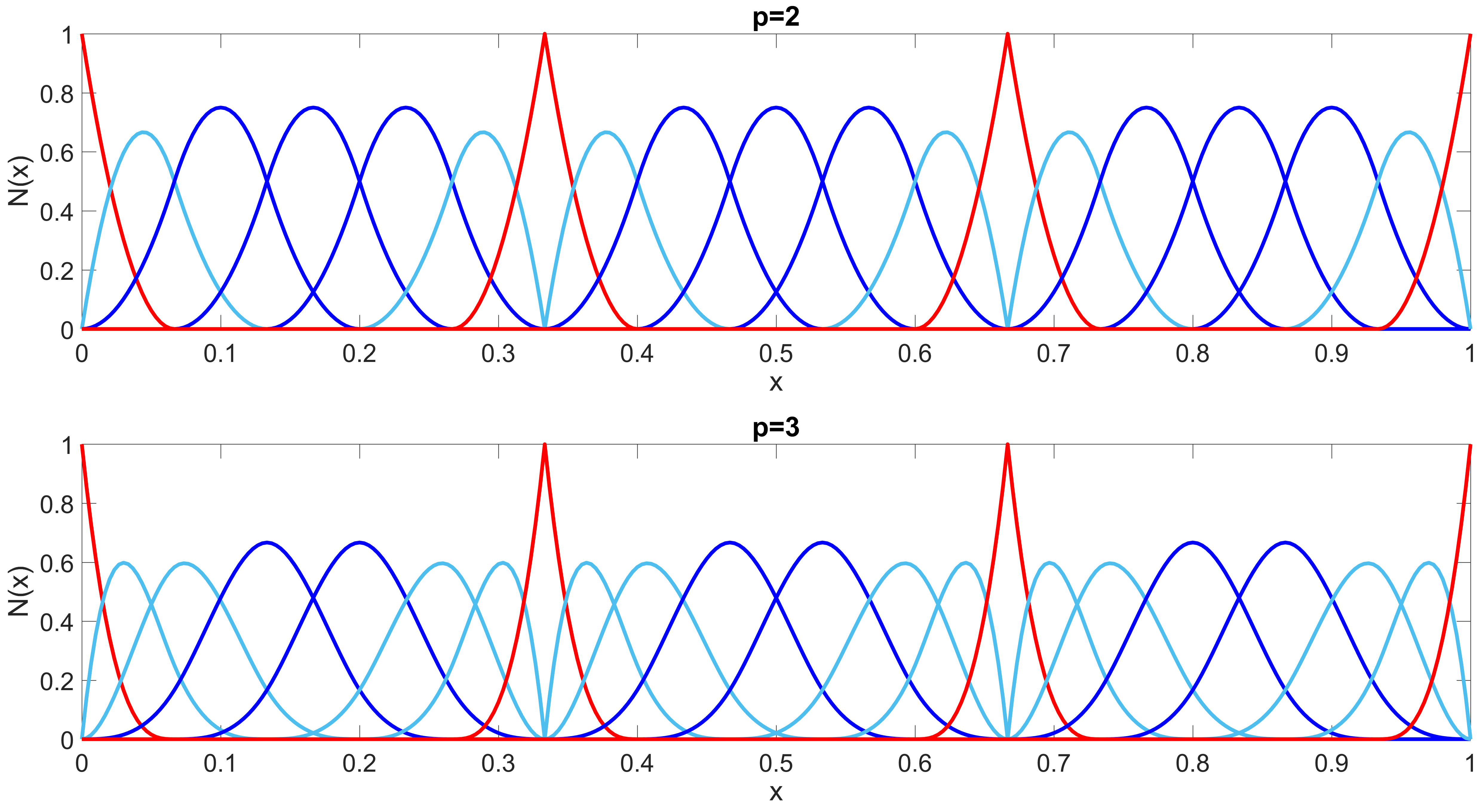}
\caption{Example of $C^2$ cubic (top) and $C^1$ quadratic (bottom) basis with two $C^0$ separators. Each block consists of homogeneous inner basis functions (blue) and $p-1$ basis functions near the boundaries of the block (light blue). Only the basis functions of the separators (red) have support over two blocks.}
\end{figure}

Figure 1 illustrates the basis of quadratic and cubic rIGA for a domain made of three blocks with five elements each interconnected by $C^0$ separators. Only the separator basis functions have support over two neighbouring blocks; all other basis functions have support only on their blocks. The knot vector contains a repetition of knots for each separator, for example, $[..., 1/3, 1/3, 1/3, ..., 2/3, 2/3, 2/3, ...]$ for the cubic case. Thus, only $C^0$-continuity is attained at these knots, while elsewhere the basis functions are $C^{p-1}$-continuous. In general, B-spline basis functions of order $p$ have $p - m_i$ continuous derivatives across knot $\xi_i$, where $m_i$ is the multiplicity of the knot. Increasing the multiplicity of existing knot values decreases the continuity and support size of the basis functions without adding new elements.

Such reduction of the inter-element continuity in rIGA improves the performance of direct solvers for a fixed logical mesh \citep{garcia2017value}. Meanwhile, if the number of elements in the original mesh is kept constant, each separator brings additional degrees of freedom to the discrete problem. As the global continuity is reduced to $C^0$, the method reduces to the standard finite-element approach with localized support of the basis functions. The increase in degrees of freedom may also increase the total solution cost of the system when compared to that of the original IGA problem. The optimal continuity reductions in terms of minimizing the number of FLOPs needed to perform the factorization of the resulting system matrix are studied in \citep{garcia2017value}. The optimal discretization decreases the cost of LU factorization while increasing the approximation properties of the discrete space.

\section{Problem formulation}
\label{S:3}

We consider an elliptic eigenvalue problem that describes the normal modes and frequencies of free structural vibration:
\begin{align} \label{eq:Original}
\begin{split}
&{-\Delta} u = \lambda u \ \ \textrm{in} \ \Omega \\
& u=0 \ \ \textrm{on} \ \partial \Omega,
\end{split}
\end{align}
where $\Delta = \nabla^2$ is the Laplace operator and $\Omega \subset \mathbb{R}^\mathit{d}$, with $d=1,2,3,$ is a bounded open domain with Lipschitz boundary. This is a classical problem in engineering, which has been thoroughly studied in the isogeometric analysis literature, e.g. \citep{cottrell2006isogeometric, cottrell2007studies, hughes2008duality}. A similar dispersion analysis has been applied to the Helmholtz equation for time-harmonic wave propagation that arises in acoustics and electromagnetics. We specialize the following derivations to the one-dimensional case to simplify the discussion. Equation \eqref{eq:Original} has an infinite set of eigenvalues ${\lambda _j} \in {\mathbb{R}^+}$ and an associated set of orthonormal eigenfunctions ${u_j}$
\begin{equation}
{0 < \lambda _1} < {\lambda _2} \le ... \le {\lambda _j} \le ...
\end{equation}
\begin{equation}
({u_j},{u_k}) = \int_{\Omega}  {{u_j}(x){u_k}} (x) \mathrm{d}x = {\delta _{jk}},
\end{equation}
where $\delta _{jk}$ is the Kronecker delta. The eigenvalues are real, positive, and countable. For each eigenvalue $\lambda _j$ there exists an eigenfunction $u_j$. The normalized eigenfunctions form an $L_2$-orthonormal basis and, as a consequence, they are also orthogonal in the energy inner product
\begin{equation}
{(\nabla{u_j},\nabla{u_k})} = ({\lambda _j}{u_j},{u_k}) = {\lambda _j}{\delta _{jk}}.
\end{equation}

The standard weak form for the eigenvalue problem is stated as follows: Find all eigenvalues ${\lambda _j} \in {\mathbb{R}^+}$ and eigenfunctions ${u_j} \in V$ such that, for all $w \in V$,
\begin{equation} \label{eq:weak}
a(w,{u_j}) = {\lambda _j}(w,{u_j}),
\end{equation}
where
\begin{equation}
a(w,{u_j}) = \int_{\Omega} \left( \nabla w \cdot \nabla u \right) \mathrm{d} x,
\end{equation}
and $V$ is a closed subspace of ${{H^1}(\Omega )}$. We use the standard notation \citep{strang1973analysis, hughes2014finite}, where $( \cdot , \cdot )$ and $a( \cdot , \cdot )$ are symmetric bilinear forms defining the following norms
\begin{equation}
\left\| w \right\|_E^2 = a(w,w),\ \ \
\left\| w \right\|_{}^2 = (w,w),
\end{equation}
for all $v,w \in V$. The energy norm is denoted as ${\left\| {\; \cdot \;} \right\|_E}$ and is equivalent to the ${{H^1_0}(\Omega )}$ norm on $V$ and $\left\| {\; \cdot \;} \right\|$ is the standard ${L_2}(\Omega )$ norm.

The Galerkin formulation of the eigenvalue problem is the discrete form of \eqref{eq:weak}: Find the discrete eigenvalues $\lambda _j^h \in {\mathbb{R}^+}$ and eigenfunctions $u_j^h \in V^h \subset V$ such that, for all ${w^h} \in {V^h} \subset V$,
\begin{equation} \label{eq:galerkin}
a({w^h},u_j^h) = \lambda _j^h({w^h},u_j^h).
\end{equation}

In matrix form, \eqref{eq:galerkin} can be written as follows
\begin{equation}
\mathbf{K} \mathbf{U} = {\lambda^h} \mathbf{M} \mathbf{U},
\end{equation}
where $\mathbf{U}$ is the matrix of eigenvectors whose \textit{j}-th column corresponds to the coefficients of the eigenfunction with respect to the basis function. The global mass and stiffness matrices $\mathbf{M}$ and $\mathbf{K}$ are defined by
\begin{equation} \label{eq:massstiff}
{M}_{i j} = \int_{\Omega} N_i(x)  N_j(x) \mathrm{d} x,\ \ \ 
{K}_{i j} = \int_{\Omega} \left[ \nabla N_i(x) \cdot \nabla N_j(x) \right] \mathrm{d}x .
\end{equation}
Here $N_i(x)$ and $N_j(x)$ are the basis functions defined as \eqref{eq:basis1}, \eqref{eq:basis2}; subscript $p$ is omitted for brevity.

\section{Error quantification}
\label{S:4}

In this section, we analyze numerically the discrete frequency spectrum and show the eigenvalue and eigenfunction errors in the $L_2$ and energy norms. The approximation of eigenvalues and eigenfunctions is fundamental for error estimates in various boundary- and initial-value problems \citep{hughes2014finite}. The errors in the discrete approximations of these problems are expressed in terms of the eigenvalue and eigenfunction errors. The Pythagorean eigenvalue error theorem \citep{strang1973analysis} states that: for each discrete mode, the eigenvalue error and the product of the exact eigenvalue and the square of the eigenfunction error in the $L_2$-norm sum to the square of the error in the energy norm, that is,
\begin{equation} \label{eq:ptorig1}
{\lambda _j^h - {\lambda _j}} + {\lambda _j}{\left\| {u_j^h - {u_j}} \right\|}^2 = {\left\| {u_j^h - {u_j}} \right\|_E^2}.
\end{equation}

In the figures below, we plot the quantities of the theorem in the following form, which shows the relative error in eigenvalues
\begin{equation} \label{eq:ptorig2}
\frac{{\lambda _j^h - {\lambda _j}}}{{{\lambda _j}}} + {\left\| {u_j^h - {u_j}} \right\|^2} = \frac{{\left\| {u_j^h - {u_j}} \right\|_E^2}}{{{\lambda _j}}}.
\end{equation}

For our dispersion studies, we consider the one-dimensional problem described in the previous section. The mesh is uniform unless otherwise specified and is chosen in such a way that the number of elements $N_e = 1000$. The number of degrees of freedom (i.e. discrete modes) in 1D isogeometric analysis with Dirichlet boundary conditions imposed at the boundaries of the domain is equal to $N_0 = N_e+p-2$. The exact eigenvalues and corresponding eigenfunctions of the 1D problem are 
\begin{equation}
{\lambda _j} = {j^2}{\pi ^2},\ \ \ {u_j} = \sqrt 2 \sin (j\pi x),
\end{equation}
for $j = {1, 2, ... }$. We sort the approximate eigenvalues $\lambda _j^h$ in ascending order and compare them to the corresponding eigenvalues of the exact operator ${\lambda _j}$. In the following figures, we show the full budget of the Pythagorean eigenvalue theorem \eqref{eq:ptorig2}, i.e. the relative eigenvalue (EV) errors $\frac{\lambda _j^h - {\lambda _j}}{\lambda _j}$, the $L_2$-norm eigenfunction (EF) errors $\left\| {{u_j^h} - u_j} \right\|^2$ and the relative energy-norm EF errors $\frac{\left\| {{u_j^h} - u_j} \right\|_E^2}{\lambda _j}$.

\begin{figure}[!ht]
\centering\includegraphics[width=1.0\linewidth]{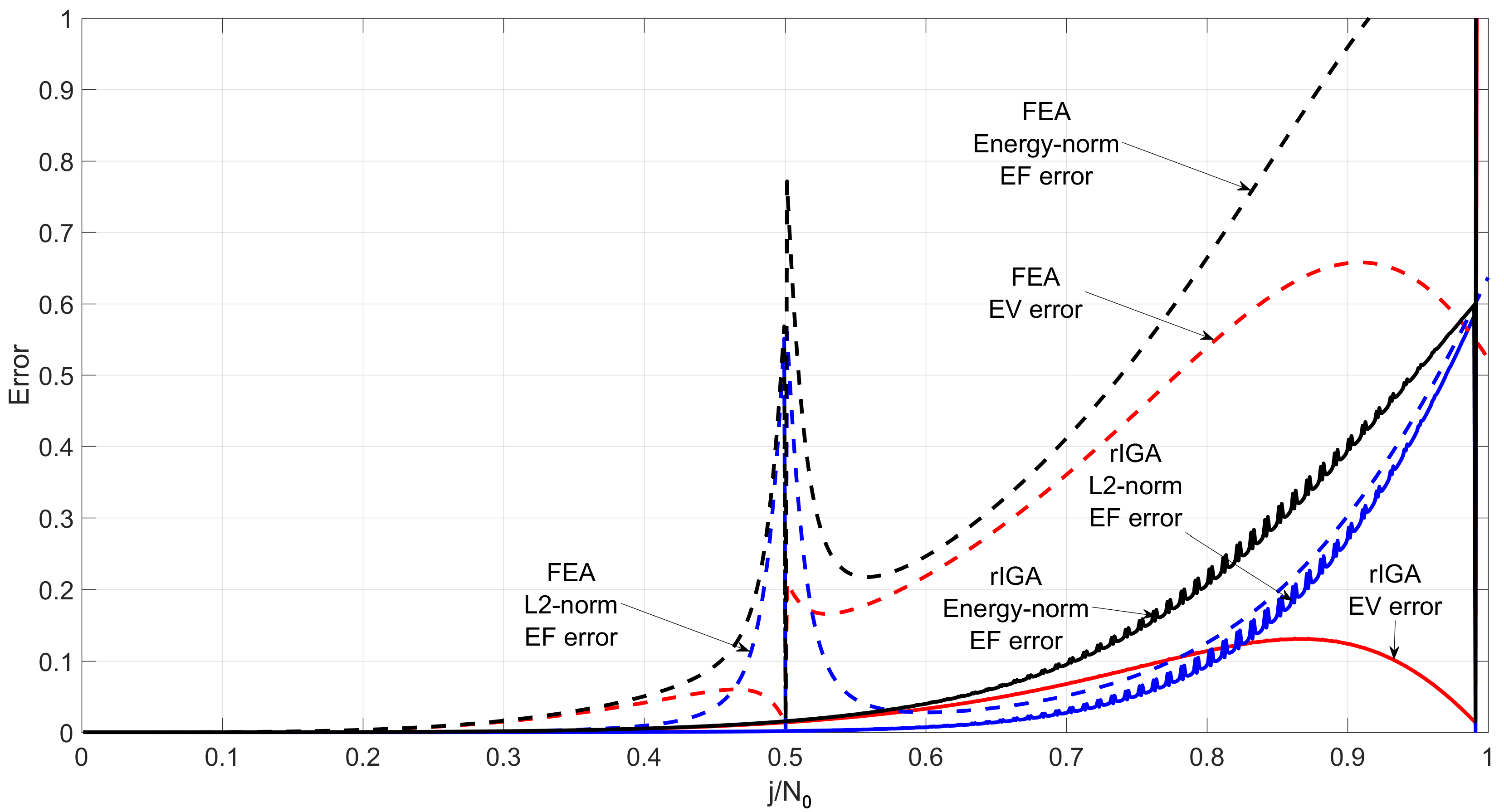}
\caption{Comparison of the eigenvalue and eigenfunction errors of quadratic $C^0$ finite elements and refined isogeometric discretization with 10 blocks of 100 $C^1$ elements each.}
\end{figure}

We first compare the accuracy of the refined isogeometric elements versus the standard finite elements. The approximate eigenvalues and eigenfunctions are significantly more accurate for IGA than for FEA for similar spatial resolutions \citep{cottrell2006isogeometric, hughes2008duality, hughes2014finite}. This improvement in the spectral accuracy of isogeometric analysis is even larger for higher-order approximations. The large spikes in the eigenfunction errors that appear at the transition points between the acoustic and optical branches of the FEA spectra, are absent in the maximum continuity discretizations. 

Figure 2 compares the approximation errors of the standard quadratic finite elements ($C^0$) with the refined isogeometric discretization that consists of 10 blocks of 100 $C^1$ elements each connected by $C^0$ separators. We assemble the mass and stiffness matrices using the standard Gaussian quadratures in this example. One can observe that the rIGA modes are more accurate compared to the FEA ones in the whole spectrum. The eigenvalue errors of rIGA are smaller than those of the FEA, expect in the outlier region at the right of the spectrum. We analyze these differences in the approximation at the high frequencies in the following sections. The convergence rate in the relative error in the eigenvalues of both methods for $h \rightarrow 0$ is the same as in IGA and FEA, that is $2p$.

For the following plots, we propose a different representation of the approximation errors that allows us to illustrate clearly how the spectrum varies when we add $C^0$ separators. The abscissa of the following error plots shows the quantity $j / N_0$, where $j$ is the number of the discrete mode. For an isogeometric discretization with uniform continuity of the basis $j / N_0 \le 1$. As the knot repetition in rIGA adds new basis functions, the number of discrete modes becomes larger than $N_0$. These new modes are shown in the right parts of the plots, where $j / N_0 > 1$. This form of representation allows us to clearly track the changes in errors for the original discrete modes and see the outlier behaviour of the new modes.

\begin{figure}[!ht]
\centering\includegraphics[width=1.0\linewidth]{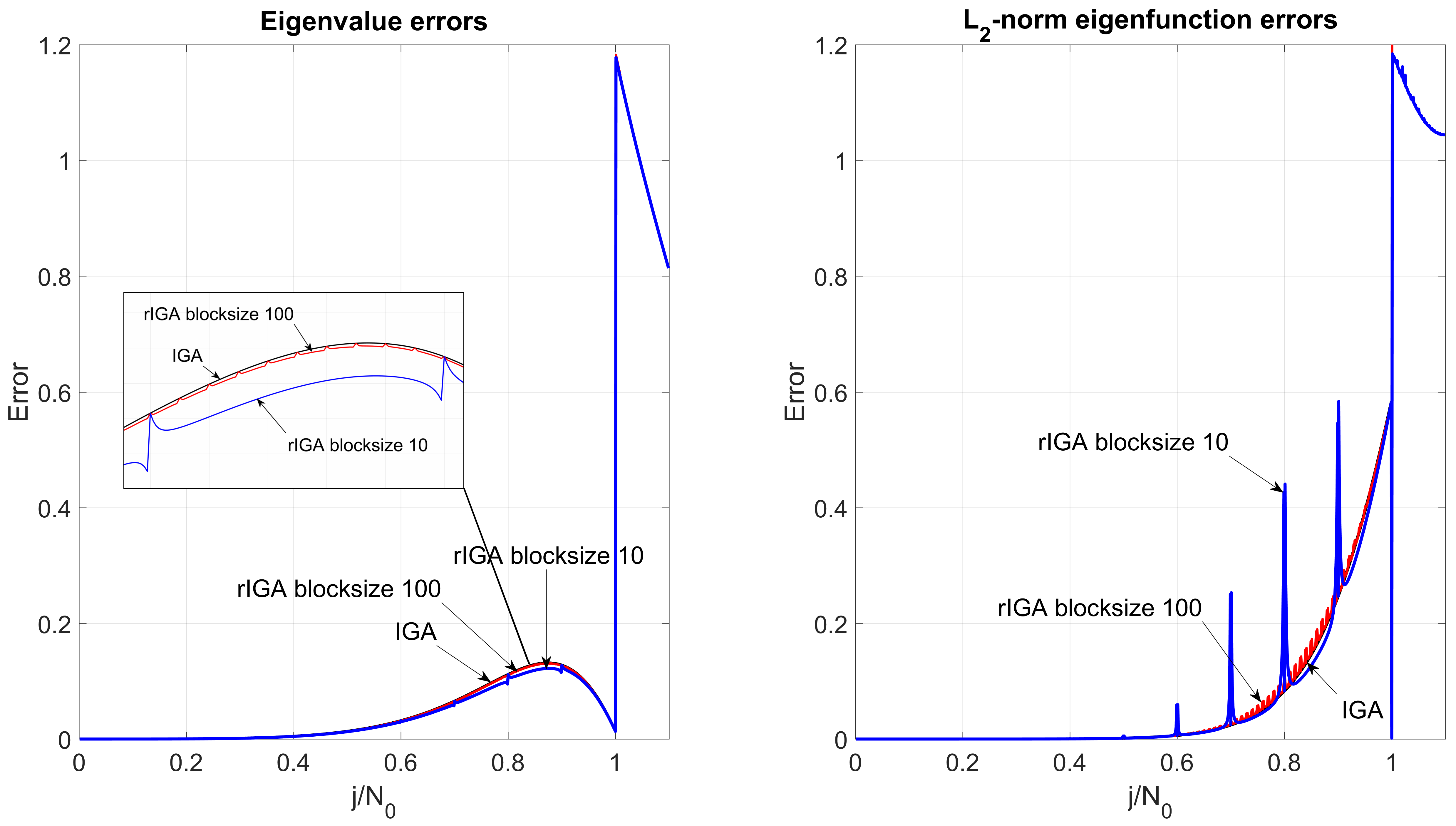}
\caption{Eigenvalue and eigenfunction errors of quadratic IGA discretizations (black line), rIGA  with 10 blocks of size 100 (red line) and rIGA with 100 blocks of size 10 (blue line).}
\end{figure}

\begin{figure}[!ht]
\centering\includegraphics[width=1.0\linewidth]{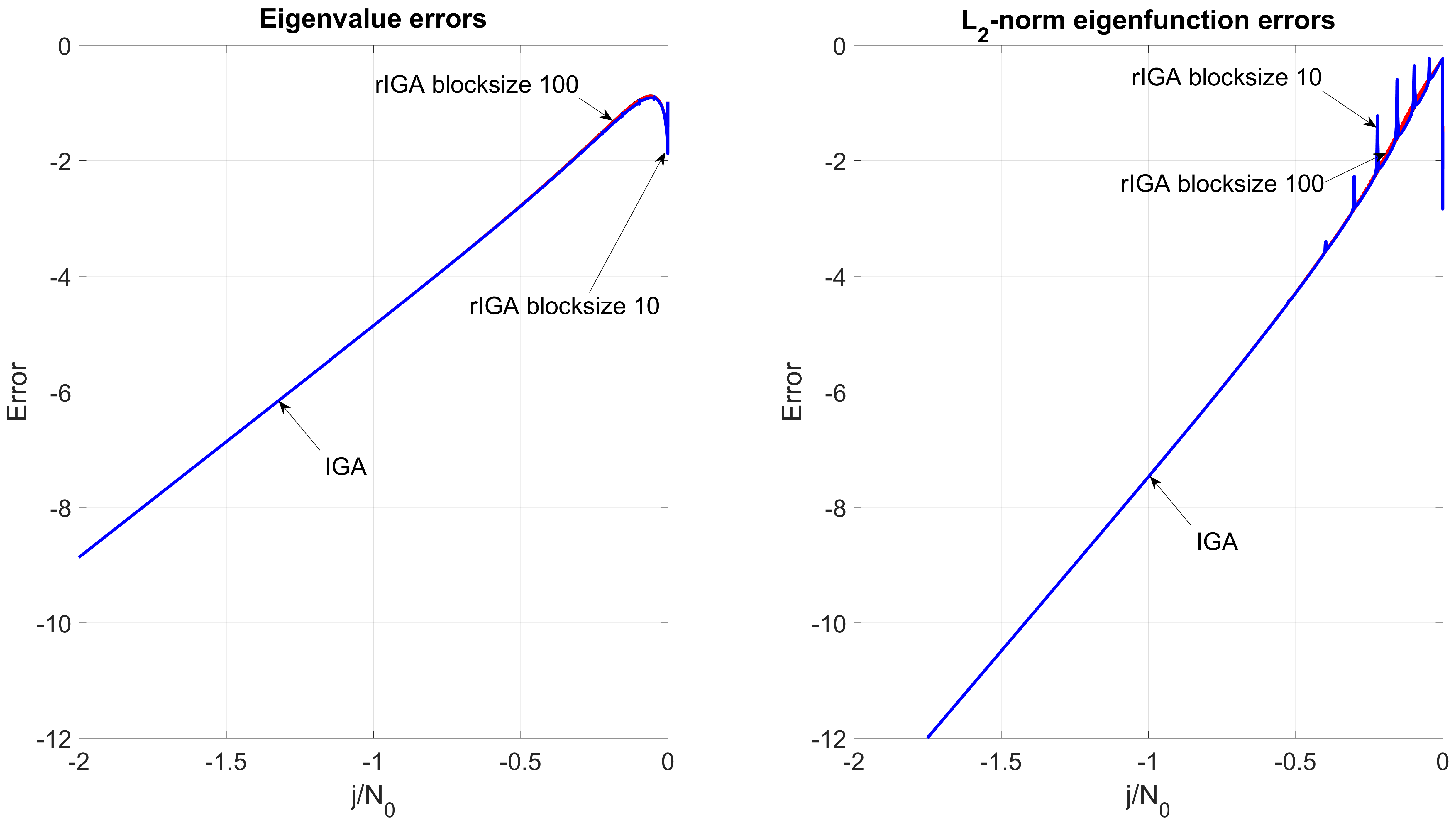}
\caption{Eigenvalue and eigenfunction errors on logarithmic scale of quadratic IGA discretizations (black line), rIGA  with 10 blocks of size 100 (red line) and rIGA with 100 blocks of size 10 (blue line).}
\end{figure}

Figures 3 and 4 show the spectra of rIGA discretizations for the same number of elements (1000), but divided into different block sizes with $C^0$ separators between the blocks. One may observe that the only difference between the spectra is the presence of spikes, which are very small in the eigenvalues and quite large in the eigenfunctions. These spikes are smaller than those in the high-order FEA spectra, but, as we show in the next section, the nature of this ``branching'' is similar. The number of branches for $j / N_0 \in [0,1]$ is equal to the size of each rIGA block (10 and 100 for the cases shown in Figures 3 and 4, respectively). The modes from the right part of the spectra ($j / N_0 > 1$) are similar to the standard outlier modes of isogeometric analysis. Their number is equal to the number of new degrees of freedom introduced by the $C^0$ separators.

For large sizes of rIGA block (e.g., 100 in Figure 3), the number of spikes in the eigenfunction errors is large, but their magnitudes are small. With the decreasing block size, the spikes gradually decrease in number, but increase in size (Figure 4). Blocks of size one produce the finite element discretization, i.e., we obtain $p$ branches of modes with large spikes between them. From the point of view of direct solvers, a very large increment in the number of separators (i.e., small $C^{p-1}$ blocks) increases the number of degrees of freedom, thus resulting in a much higher factorization cost, especially for 3D problems \citep{garcia2017value}. Various configurations are tested in \citep{garcia2017value} considering both the number of FLOPs required to eliminate the degrees of freedom and the actual computational time to factorize the matrix. Block sizes of several dozens of degrees of freedom tend to deliver the best performance, although the regions with minimal computational cost are fairly wide, varying between 20 and 150 elements per block. The optimum size tends to be around 20 elements per block.

Further increase in the number of $C^0$ separators between the blocks of the $C^{p-1}$ continuity increases the branching of the rIGA spectrum. Figure 5 shows the eigenvalue errors of quadratic rIGA with the block sizes of eight, four, two, and one (the latter corresponds to the $C^0$ FEA discretization). The number of branches in the left part of the spectrum ($j / N_0 \in [0,1]$) is equal to the rIGA block size. As the block size decreases, these errors converge to the acoustic branch of the FEA spectrum. In the right part of the spectrum ($j / N_0 > [0,1]$), we observe some similarity between the rIGA eigenvalue errors for the block size of two and the FEA optical branch.

\begin{figure}[!ht]
\centering\includegraphics[width=1.0\linewidth]{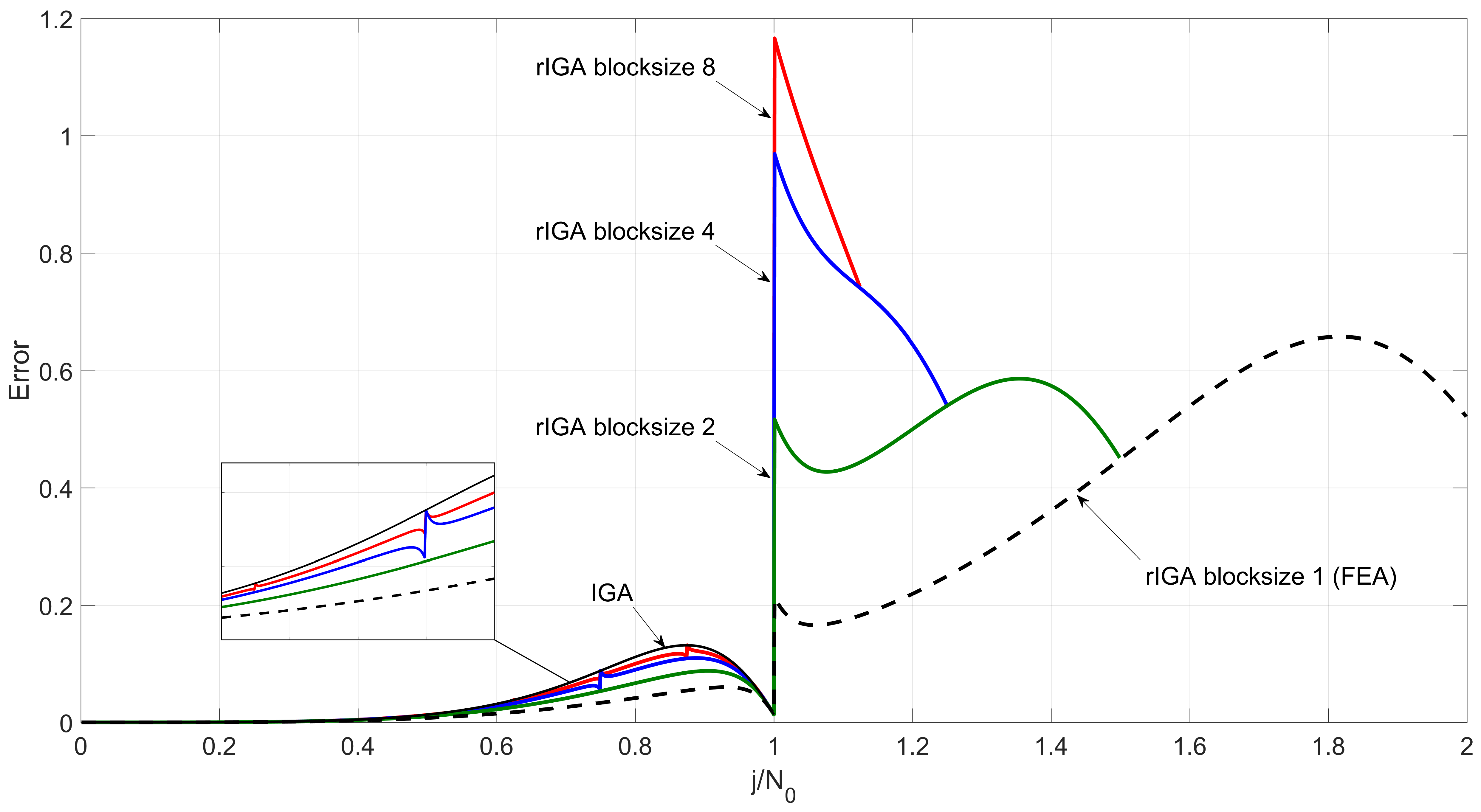}
\caption{Eigenvalue errors of quadratic rIGA discretizations with the block sizes of eight (blue line), four (red line), two (green line), and one (dashed black line). IGA eigenvalue errors are shown for comparison (solid black line).}
\end{figure}

\section{Stopping bands in eigenvalue errors}
\label{S:5}

The heterogeneity of the high-order finite element basis functions leads to branching of the discrete spectrum and a fast degradation of the accuracy for higher frequencies (optical branches). B-spline $C^{p-1}$ basis functions are homogeneous on uniform meshes and do not exhibit such branching patterns (with an exception of the outliers that correspond to the basis functions with support on the boundaries of the domain for higher $p$).

Stopping bands in eigenvalue errors for several high-order finite element methods were studied by \citet{thompson1994complex}. They extended the standard dispersion analysis technique to include complex wavenumbers and showed that some solutions are not purely propagating (real wavenumbers) but are attenuated (complex wavenumbers) waves. The stiffness and mass matrices of high-order finite elements can be partitioned into the matrix block form with one partition corresponding to interactions among bubble (internal to the element) shape functions. The elimination of variables in equations corresponding to the interior degrees of freedom is performed by means of Schur complements \citep{thompson1994complex}.

In this section, we show that the nature of the stopping bands in FEA and rIGA is similar and describe these stopping bands using the Schur complements without recurring to the introduction of complex wave numbers, which seems unnecessary. Indeed, the structure of the global mass and stiffness matrices for rIGA with many blocks connected the their neighbours by $C^0$ separators is similar to the case of high-order finite elements \citep{garcia2017value, garcia2017optimally}. Figure 6 compares the sparsity patterns of the mass matrices for $C^0$ elements of 11th order and quadratic rIGA with ten elements in the inner blocks. In both cases, we assemble blocks with ten degrees of freedom which are internal to the $C^0$ separators and thus the resulting spectra have ten spikes. The high-order finite element method produces a denser matrix than the lower-order rIGA.

\begin{figure}[!ht]
\centering\includegraphics[width=1.0\linewidth]{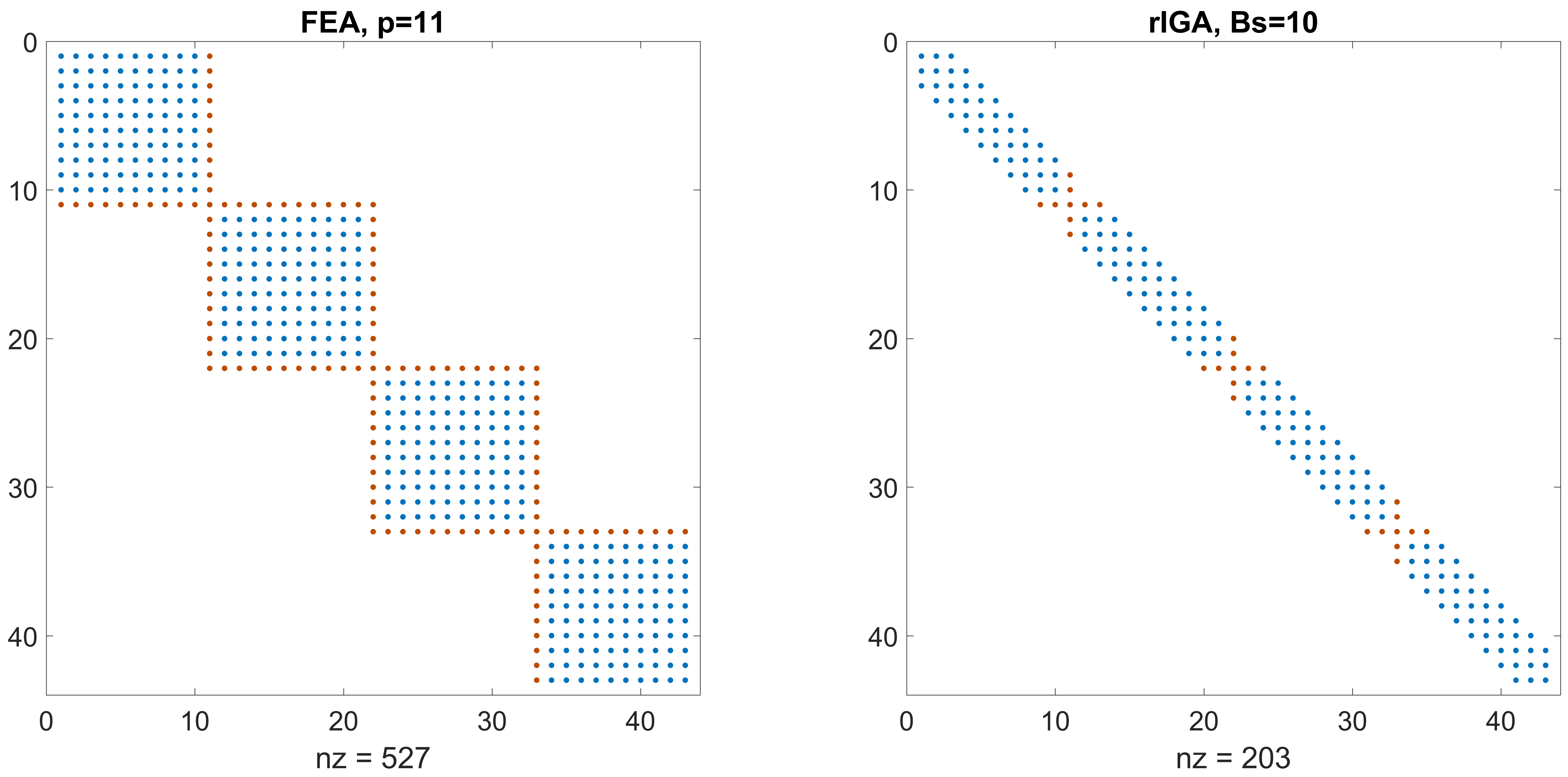}
\caption{Matrix sparsity patterns for $C^0$ elements of 11th order (left) and quadratic $C^1$ elements with $C^0$ separators between the blocks of 10 elements (right). The separator degrees of freedom are highlighted in red.}
\end{figure}

We write the original matrix equation for the eigenvalue problem as follows

\begin{equation} \label{eq:start}
\left[
\left( \begin{array}{cc}
\mathbf{K}_{bb} & \mathbf{K}_{bi} \\
\mathbf{K}_{ib} & \mathbf{K}_{ii} \end{array} \right)
- \lambda^h
\left( \begin{array}{cc}
\mathbf{M}_{bb} & \mathbf{M}_{bi} \\
\mathbf{M}_{ib} & \mathbf{M}_{ii} \end{array} \right)
\right]
\left( \begin{array}{cc}
\mathbf{U}_b \\
\mathbf{U}_i \end{array} \right)
= \mathbf{0}.
\end{equation}
Here ``b'' stands for ``bubble'' and ``i'' stands for ``interface'' basis functions (see also Figure 1). $\mathbf{K}_{bb}$ and $\mathbf{M}_{bb}$ consist of the bubble degrees of freedom of the stiffness and mass matrices, respectively. We rewrite \eqref{eq:start} as
\begin{align} \label{eq:system2}
\begin{split}
\begin{cases}
\left( \mathbf{K}_{bb} - \lambda^h \mathbf{M}_{bb} \right) \mathbf{U}_b + \left( \mathbf{K}_{bi} - \lambda^h \mathbf{M}_{bi} \right) \mathbf{U}_i &= \mathbf{0},\\
\left( \mathbf{K}_{ib} - \lambda^h \mathbf{M}_{ib} \right) \mathbf{U}_b + \left( \mathbf{K}_{ii} - \lambda^h \mathbf{M}_{ii} \right) \mathbf{U}_i &= \mathbf{0}.
\end{cases}
\end{split}
\end{align}

\citet{thompson1994complex} expressed the bubble degrees of freedom in terms of the interfacial ones. From the first equation of \eqref{eq:system2} we write
\begin{equation}
\mathbf{U}_b = -(\mathbf{K}_{bb}  - \lambda^h \mathbf{M}_{bb})^{-1} (\mathbf{K}_{bi} - \lambda^h \mathbf{M}_{bi}) \mathbf{U}_i,
\end{equation}
to eliminate $\mathbf{U}_b$ from \eqref{eq:system2} and express it in terms of $\mathbf{U}_i$. This process is appropriate as long as $\det(\mathbf{K}_{bb}  - \lambda^h \mathbf{M}_{bb})\neq \mathbf{0}$. However, this condition is not satisfied at the stopping bands. At the stopping bands we have

\begin{equation}
\det(\mathbf{K}_{bb}  - \lambda^h_b \mathbf{M}_{bb}) = \mathbf{0},
\end{equation}
for particular eigenvalues $\lambda^h_b$, which are the eigenvalues of the local ``bubble'' subsystems, as well as the eigenvalues of the global system. 

Therefore, 
\begin{equation}
\left( \mathbf{K}_{bb}-\lambda^h_b \mathbf{M}_{bb}\right)\mathbf{U}_b=\mathbf{0}.
\end{equation}
Additionally, the local bubble eigenvalue problem can be expressed as 
\begin{equation}
\left( \mathbf{K}^L_{bb}-\lambda^h_b \mathbf{M}^L_{bb}\right)\mathbf{U}^L_b=\mathbf{0},
\end{equation}
where the superscript $L$ denotes the restriction of the mass and stiffness matrices as well as the degrees of freedom to the local block problem. In FEA, the local problem block corresponds to the interior bubbles of each individual higher-order finite element. In rIGA, the local problem refers to all the interior bubbles that are zero at the lower continuous separator. In this case, $\det(\mathbf{K}_{ii}  - \lambda^h_b \mathbf{M}_{ii})\neq \mathbf{0}$ and from the second equation of \eqref{eq:system2} we can express the interface degrees of freedom as
\begin{equation} \label{eq:system3}
\mathbf{U}_i = -(\mathbf{K}_{ii}  - \lambda^h_b \mathbf{M}_{ii})^{-1} (\mathbf{K}_{ib} - \lambda^h_b \mathbf{M}_{ib}) \mathbf{U}_b.
\end{equation}

By substituting this expression into \eqref{eq:system2}, we obtain
\begin{equation} \label{eq:system4}
(\mathbf{K}_{bb}  - \lambda^h_b \mathbf{M}_{bb}) \mathbf{U}_b - ( \mathbf{K}_{bi} - \lambda^h_b \mathbf{M}_{bi} ) (\mathbf{K}_{ii}  - \lambda^h_b \mathbf{M}_{ii})^{-1} (\mathbf{K}_{ib} - \lambda^h_b \mathbf{M}_{ib}) \mathbf{U}_b = \mathbf{0}.
\end{equation}

The first term of \eqref{eq:system4} is zero at the stopping bands as explained above. Thus, we write \eqref{eq:system4} as
\begin{equation} \label{eq:system5}
( \mathbf{K}_{bi} - \lambda^h_b \mathbf{M}_{bi} ) (\mathbf{K}_{ii}  - \lambda^h_b \mathbf{M}_{ii})^{-1} (\mathbf{K}_{ib} - \lambda^h_b \mathbf{M}_{ib}) \sum_{L} {u}_{b}^L \alpha^L = \mathbf{0},
\end{equation}
where ${u}_{b}^e$ are the local (to the element for FEA and block for rIGA) bubble eigenfunctions that correspond to each $\lambda^h_b$ at local problem such that $\mathbf{U}_b = \sum_{e} {u}_{b}^e \alpha^e$. The local eigenvalue problem defines the shape of each bubble eigenfunction. Thus, we use \eqref{eq:system5} to set the orientation of each local bubble to form a global wave. We use Galerkin's projection to set up and solve the reduced system to determine the values of $\alpha^L$. We assemble $\mathbf{U}_b$ from the expression above using the local bubbles and $\alpha^L$, and then use \eqref{eq:system3} to compute the remaining interface degrees of freedom. The ``bubble'' matrices $\mathbf{K}_{bb}$ and $\mathbf{M}_{bb}$ define the stopping bands. A high-order finite element spectrum has $p-1$ stopping bands, while rIGA spectrum has $Bsize + p - 2$ stopping bands, where $Bsize$ is the size of a $C^{p-1}$ block.

\section{Outliers}
\label{S:6}

In this section, we unify the description of the outlier modes in IGA and rIGA discretizations. Outliers are the modes with large errors that appear at the high frequency part of the spectrum in isogeometric analysis. This phenomenon is related to the types of basis functions and the discrete equations they produce. In IGA for a uniform open knot vector, all basis functions are spatial translations of one another except those associated with the open knot vectors. These clustered knots are responsible for the outliers \citep{cottrell2006isogeometric}. As shown in Section 4, rIGA discretizations have outlier modes with large errors at the high end of the spectrum ($j / N_0 > 1$). Let $N_{sep}$ be the number of $C^0$ separators in a 1D rIGA discretization. Table 1 lists the number of outlier frequencies present in the maximum continuity IGA and the variable continuity rIGA spectra for different numbers of separators and polynomial orders, while we assume that the interior basis in each block has maximum continuity.

\begin{center}
\begin{tabular}{ |c|c|c| } 
 \hline
 \rowcolor{lightgray} \textit{p} & IGA & rIGA \\
 \hline
 \textit{2} & 0 & 0+$N_{sep}$ \\ 
 \hline
 \textit{3} & 2 & 2+$2 N_{sep}$ \\ 
 \hline
 \textit{4} & 2 & 2+$3 N_{sep}$ \\ 
 \hline
 \textit{5} & 4 & 4+$4 N_{sep}$ \\ 
 \hline
 \textit{6} & 4 & 4+$5 N_{sep}$ \\ 
 \hline
 \textit{7} & 6 & 6+$6 N_{sep}$ \\ 
 \hline
 \textit{8} & 6 & 6+$7 N_{sep}$ \\ 
 \hline
\end{tabular}
\captionof{table}{Number of outliers in $C^{p-1}$ IGA and rIGA for different number of separators and polynomial order.}
\end{center}

The total number of outliers in both IGA and rIGA is
\begin{equation} \label{eq:outliernum}
N_{out} = N_{out}^{IGA}  + (p-1) N_{sep},
\end{equation}
where $N_{out}^{IGA}$ is the number of outliers in $C^{p-1}$ isogeometric analysis. This number is different for Dirichlet and Neumann boundary conditions, though their pattern of appearance is similar. There are $p-1$ or $p$ basis functions near each boundary of the domain for, respectively, Dirichlet and Neumann boundary conditions. These basis functions contribute to the outliers, however, not all of them produce an outlier. For Dirichlet boundary conditions, $C^{p-1}$ isogeometric elements produce two outliers for each odd polynomial order $p$ starting from cubics. For Neumann boundary conditions, we have two new outliers for each even polynomial order $p$. For example, quadratic isogeometric elements have no outliers for the Dirichlet case and two outliers for the Neumann case; cubics have two outliers for both types of boundary conditions.

\begin{figure}[!ht]
\centering\includegraphics[width=1.0\linewidth]{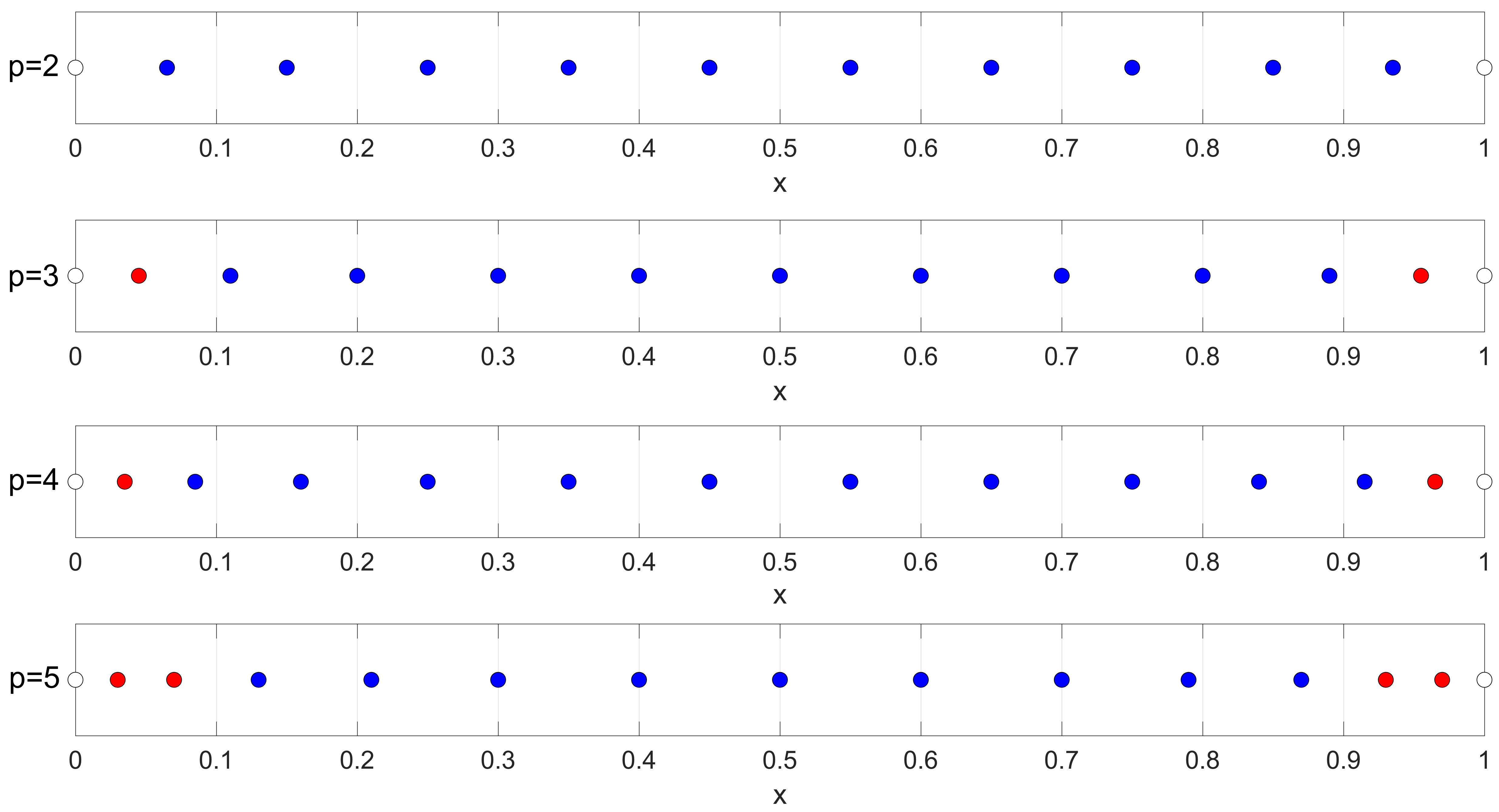}
\caption{Distribution of control points in quadratic IGA. Blue: inner points; red: extra points that cause outliers; white: boundary points.}
\end{figure}

We relate the outliers in isogeometric analysis to the knot distribution. Figure 7 shows the distribution of control points in IGA of higher continuity for different order of elements. There are extra points that create the denser regions near the boundaries. All blue points near the boundaries can be shifted such that they form a uniform knot spacing, but there is no space for red extra points. This may illustrate why the number of outliers is different for odd and even $p$. Each odd $p$ brings a knot from each side that does not fall into the previous distribution of knots. In \citep{calo2017quadrature} the authors show that the IGA spectrum on special non-uniform meshes (made of two types of elements, one of which is much larger than the other) is similar to the finite-element spectrum for a uniform mesh. This is not surprising since the basis functions for this IGA discretization are similar to the $C^0$ basis on a uniform mesh. Thus, the knot spacing controls the shape of the basis functions which, in turn, affects the spectrum.

\begin{figure}[!ht]
\centering\includegraphics[width=1.0\linewidth]{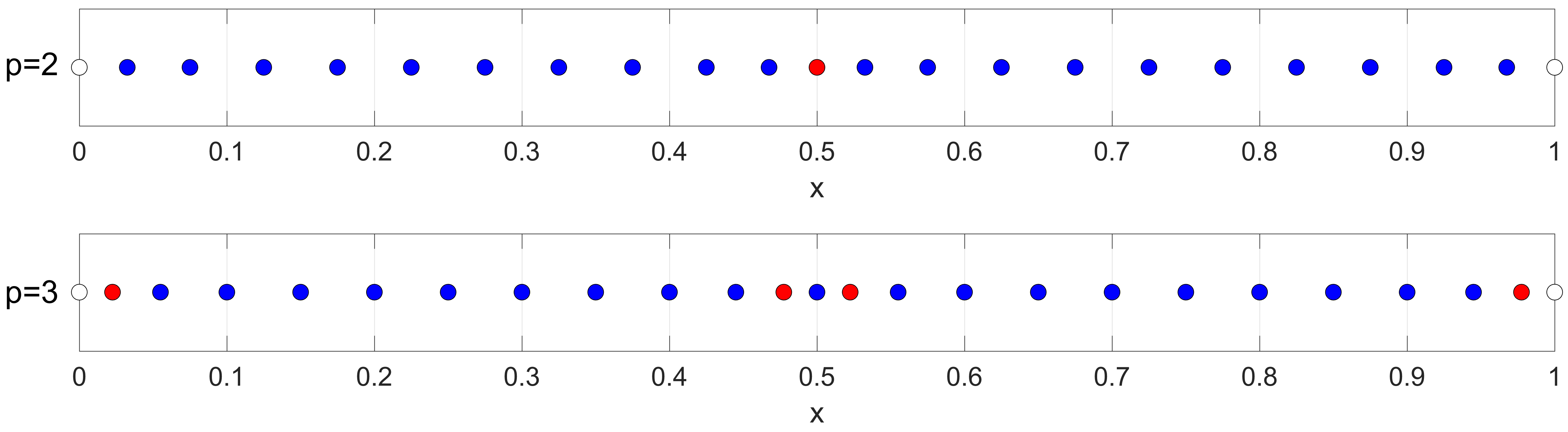}
\caption{Distribution of control points in quadratic rIGA with two blocks of ten elements. Blue: inner points; red: extra points that cause outliers; white: boundary points.}
\end{figure}

Figure 8 shows the distribution of control points in rIGA with one separator at the center. For quadratic elements, the added point does not fit on a uniform distribution of the control points and creates one outlier. For cubic elements, we add two extra points in the center and hence they introduce two outliers. The same happens for any higher-order rIGA. The number of outliers is equal to the number of these extra control points (and their associated basis functions).

Let us now consider the shape and accuracy of the eigenfunctions in rIGA. \citet{hughes2014finite} studied the outlier eigenfunctions in high-order IGA, which are zero throughout most of the domain, and only nonzero near the domain boundaries. Figure 9 shows five eigenfunctions of the quadratic rIGA with two $C^0$ separators between blocks of 64 $C^1$ elements. The last two eigenfunctions correspond to the outliers in this case (according to \eqref{eq:outliernum}, the number of outliers is equal to the number of separators for quadratic elements). Frequency analysis shows that, while most eigenfunctions have a dominant frequency (that corresponds closely to the frequency of the exact solution), the outlier eigenfunctions contain many frequencies without clear peaks.

\begin{figure}[!ht]
\centering\includegraphics[width=1.0\linewidth]{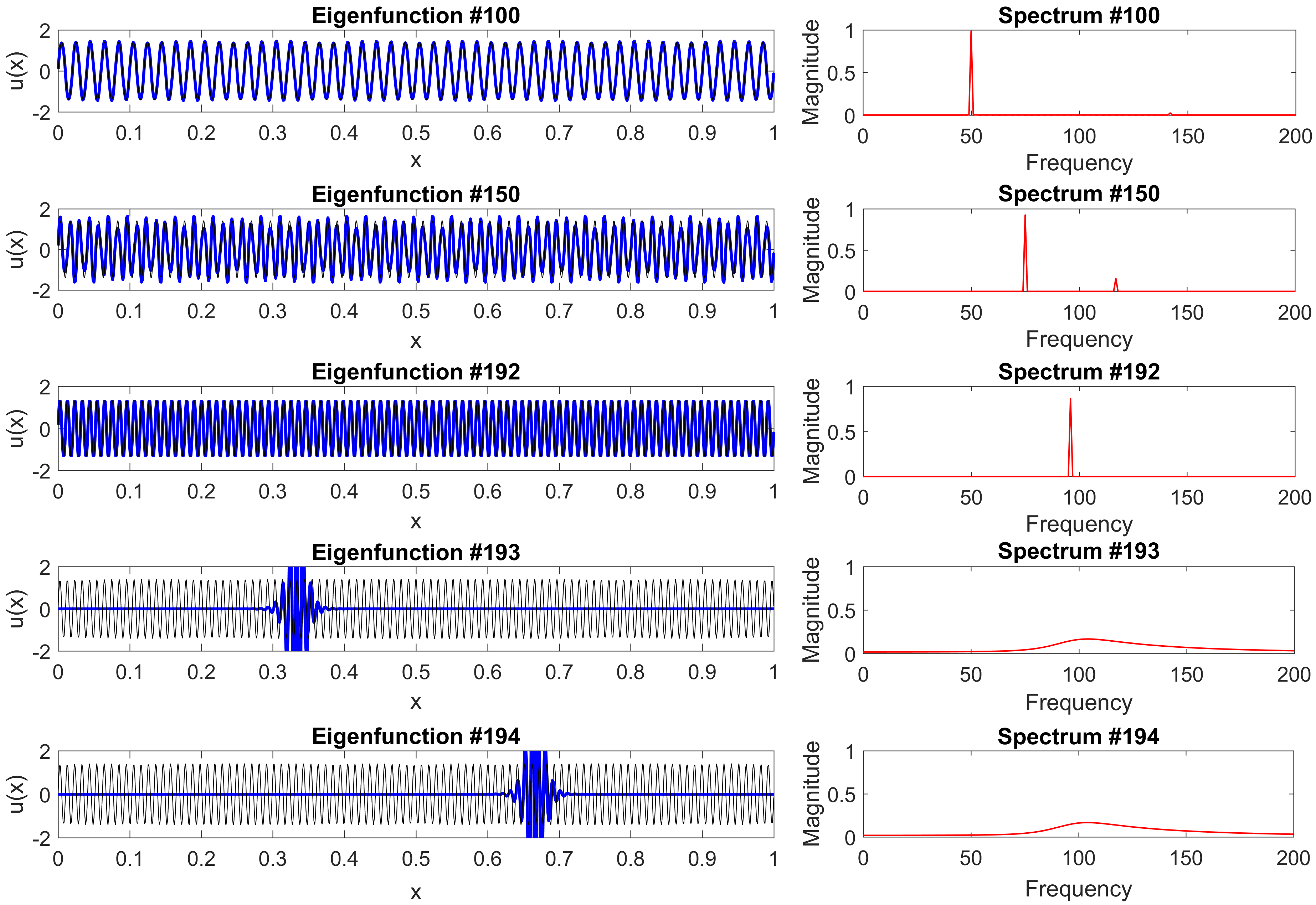}
\caption{The 100th, 150th, 192nd, 193rd, and 194th exact (black) and approximate (blue) eigenfunctions of the quadratic rIGA with two separators. The last two eigenfunctions correspond to the outliers. The right panels show the frequency content of the corresponding eigenfunctions.}
\end{figure}

Figure 10 shows the six highest eigenfunctions of the cubic rIGA with one separator at the center of the domain. The first and the third outliers (192nd and 194th eigenfunctions, respectively) are the standard IGA outliers which are nonzero only near the domain boundaries. The other two (193rd and 195th eigenfunctions) are the rIGA outliers ($p-1$ per separator) which are nonzero only near the separator. The shape and frequency content of the outliers is very similar since all of them are caused by the same reason -- non-uniformity of the knot pattern. The outlier eigenfunctions, similar to the previous case, contain many frequencies and are completely spurious solutions.

\begin{figure}[!ht]
\centering\includegraphics[width=1.0\linewidth]{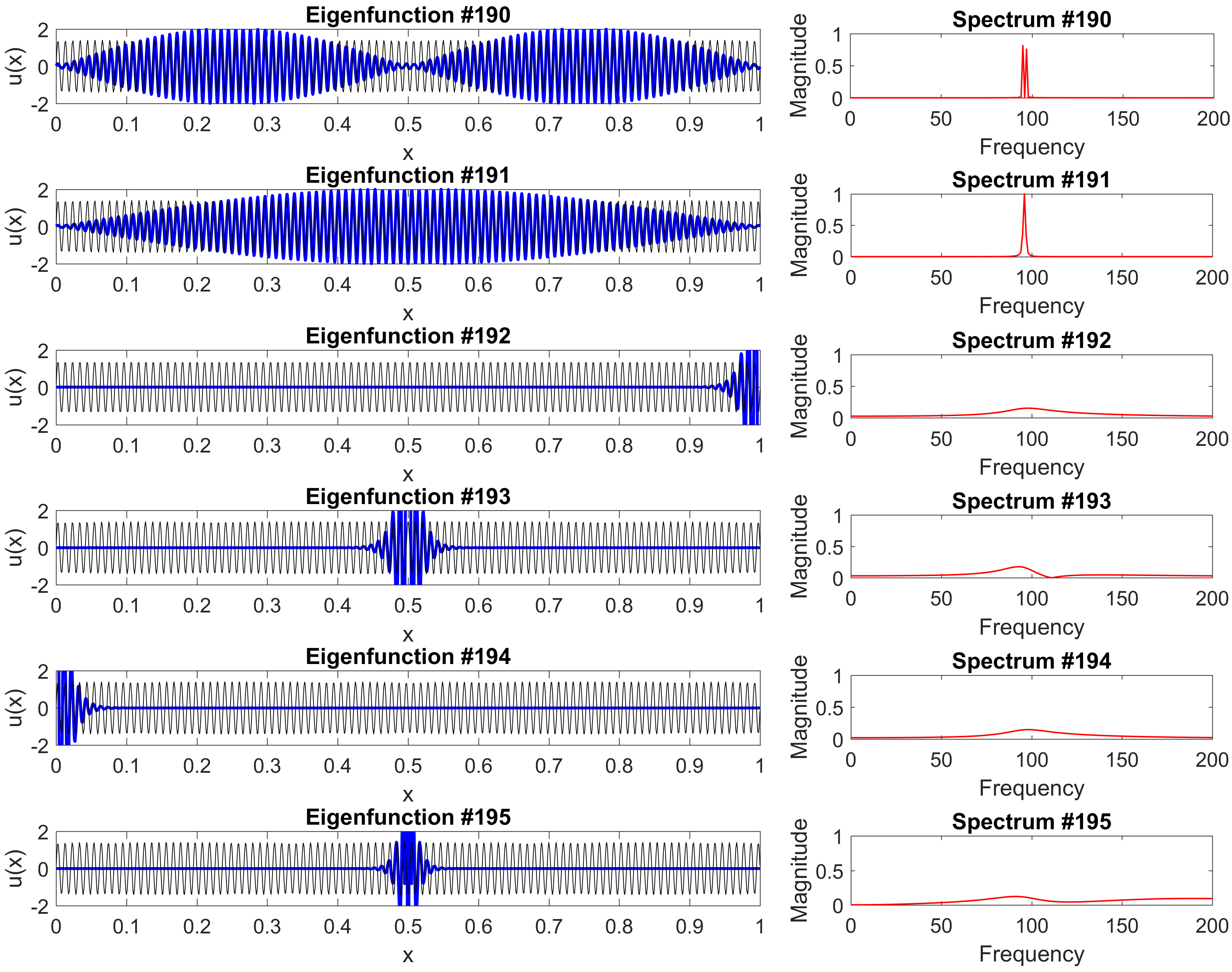}
\caption{The six highest exact (black) and approximate (blue) eigenfunctions of the cubic rIGA with one separator. The last four eigenfunctions correspond to the outliers. The right panels show the frequency content of the corresponding eigenfunctions.}
\end{figure}

The shape of the highest non-outlier eigenfunctions (e.g., \#190 and \#191 in Figure 10) is similar to amplitude-modulated (AM) waves in radio. These waves can be represented as a sum of several sine and cosine waves. As can be seen in the frequency spectra, approximate eigenfunctions consist mostly of two harmonic waves. Isogeometric eigenfunctions of even degrees can be approximately represented via sines \citep{thompson1994complex}

\begin{equation} \label{eq:u1}
u(x) = A_1 \sin({2 \pi f_1 x}) - A_2 \sin({2 \pi f_2 x}),
\end{equation}
where amplitudes $A_1$ and $A_2$ are the y-axis values from the frequency spectra in Figure 9, and frequencies $f_1$ and $f_2$ are their x-axis values. The frequencies are linked with the number of elements in the mesh as $f_2 = n - f_1$, where $n$ is the number of degrees of freedom. For even-degree isogeometric elements, the last eigenfunction (\#192 in Figure 9) is represented very well. Its frequency spectrum has only one main frequency that is the exact eigenfrequency.

Isogeometric eigenfunctions of odd degrees can be approximately represented in a similar way via cosines 

\begin{equation} \label{eq:u2}
u(x) = A_1 \cos({2 \pi f_1 x}) + A_2 \cos({2 \pi f_2 x}).
\end{equation}
For odd-degree elements, all eigenfunctions starting from the 192nd are outliers (Figure 10). For Neumann boundary conditions, the sine and cosine functions in \eqref{eq:u1} and \eqref{eq:u2} swap.

The outliers are ``extra eigenfunctions'' that do not ``fit'' the uniformly sampled spectrum. Since the number of elements in both examples shown in Figures 9 and 10 is 192, the highest wave that can propagate through this mesh has 96 full periods. The higher frequencies that the new knots add ($j/N_0 > 1$), the Nyquist sampling rate is not satisfied. That is, for lossless representation, the sampling rate must be at least twice the maximum frequency. High-order isogeometric analysis adds extra knots at the ends of the domain, but not to the internal elements. Thus, IGA leads to outlier modes that are absent in standard finite elements on uniform meshes where nodes are added to each element and the Nyquist sampling rate is not broken. The number of propagating waves is limited and the waves that do not fit into this group are undersampled waves which are identical to the outlier modes in the isogeometric spectrum. These modes behave like evanescent waves as Figures 9 and 10 show. The outliers have large peaks in the regions with the dense knot distribution (i.e., at the boundaries and separators) and are almost zero in the rest of the domain.

\section{rIGA with optimal blending}
\label{S:7}

In this section, we discuss the use of dispersion-minimizing quadratures to improve the accuracy of the rIGA approximations. The integrals in the mass and stiffness matrices \eqref{eq:massstiff} are evaluated with the help of quadrature rules. The classical IGA and FEA typically employ a Gauss quadrature that fully integrates \eqref{eq:massstiff}. Reduced integration can be used to improve the accuracy and efficiency of a numerical method. In particular, reduced integration may minimize the numerical dispersion of the finite element and isogeometric analysis approximations. Higher accuracy can be obtained using optimal blending schemes (blending the Gauss and Lobatto quadratures \citep{ainsworth2010optimally, puzyrev2017dispersion, deng2017quadratures}) or dispersion-minimizing quadrature rules, which can be constructed to produce equivalent results \citep{deng2017quadratures2, barton2017bookchapter}. This dispersion-minimizing integration improves the convergence rate of the resulting eigenvalues by two orders when compared against the fully-integrated finite, spectral, or isogeometric elements, while preserving the optimal convergence of the eigenfunctions.

To show how this improvement in the convergence rate is achieved, we consider as an example the approximate eigenvalues written as a series in $h$ for the quadratic isogeometric elements using, respectively, the Gauss and Lobatto quadrature rules with 3 quadrature points each
\begin{equation} \label{eq:iga21}
 \frac{\lambda^h_{G}}{\lambda} = 1 - \frac{1}{12} \frac{1}{5!} \lambda^2 h^4 + O(h^6),
\end{equation}
\begin{equation} \label{eq:iga22}
 \frac{\lambda^h_{L}}{\lambda} = 1 + \frac{1}{24} \frac{1}{5!} \lambda^2 h^4 + O(h^6).
\end{equation}
We blend these two schemes using a blending parameter $\tau$ leading to the approximate eigenvalues of the following form
\begin{equation} \label{eq:blend}
 \frac{\lambda^h_{B}}{\lambda} = 1 + (3 \tau - 2) \frac{1}{24} \frac{1}{5!} \lambda^2 h^4 + O(h^6).
\end{equation}

The choice of $\tau = 2/3$ allows the second term of the right-hand side of \eqref{eq:blend} to vanish and increases by two additional orders of accuracy the eigenvalue approximation when compared with the standard method. Similarly, we can improve the eigenvalue approximation quality of high-order schemes by two orders of accuracy by removing the leading order term from the error expansion. For more details, we refer the reader to \cite{puzyrev2017dispersion, deng2017quadratures}; a similar technique was used for the finite and spectral $C^0$ elements in \citep{ainsworth2010optimally}. The convergence rate for the eigenvalue errors of the optimally-blended IGA approximations is $O\left( \Omega^{2p+2} \right) $ versus the standard IGA that has a $O\left( \Omega^{2p} \right) $ convergence rate.

The standard Pythagorean eigenvalue error theorem \eqref{eq:ptorig2} assumes that the discrete method fully reproduces the inner product. To quantify the approximation errors in the case when the integrals in the Galerkin formulation are underintegrated (i.e., modified discrete inner product representations), we use the Pythagorean eigenvalue error theorem for these modified inner products \citep{puzyrev2017dispersion}:
\begin{equation} \label{eq:ptfinal}
\frac{\left\| {{u_j} - v_j^h} \right\|_E^2}{\lambda _j} = \frac{\mu _j^h - {\lambda _j}}{\lambda _j} + \left\| {{u_j} - v_j^h} \right\|^2 + \frac{\left\| {v_j^h} \right\|_E^2 - \left\| {v_j^h} \right\|_{E,h}^2}{\lambda _j} + \left(1 - \left\| {v_j^h} \right\|^2 \right),
\end{equation}
where $\mu_j^h$ and $v_j^h$ are, respectively, the discrete eigenvalues and eigenfunctions resulting from modified inner-product discretizations. The third term of \eqref{eq:ptfinal} is the error in the discrete energy norm, which is zero in the 1D case to preserve the optimal convergence in the energy norm. The last term is the error in the $L_2$ inner product, which is not zero as shown on the following figures. When the inner products are fully integrated, \eqref{eq:ptfinal} naturally reduces to the standard Pythagorean eigenvalue error theorem.

We now show the performance of these optimal quadratures in rIGA, which reduce the errors in the approximation of the eigenvalues (and, in some cases, the eigenfunctions). Similarly to optimal FEA and IGA, optimal rIGA has two extra orders of convergence when compared to the fully-integrated case for a given polynomial order.

\begin{figure}[!ht]
\centering\includegraphics[width=1.0\linewidth]{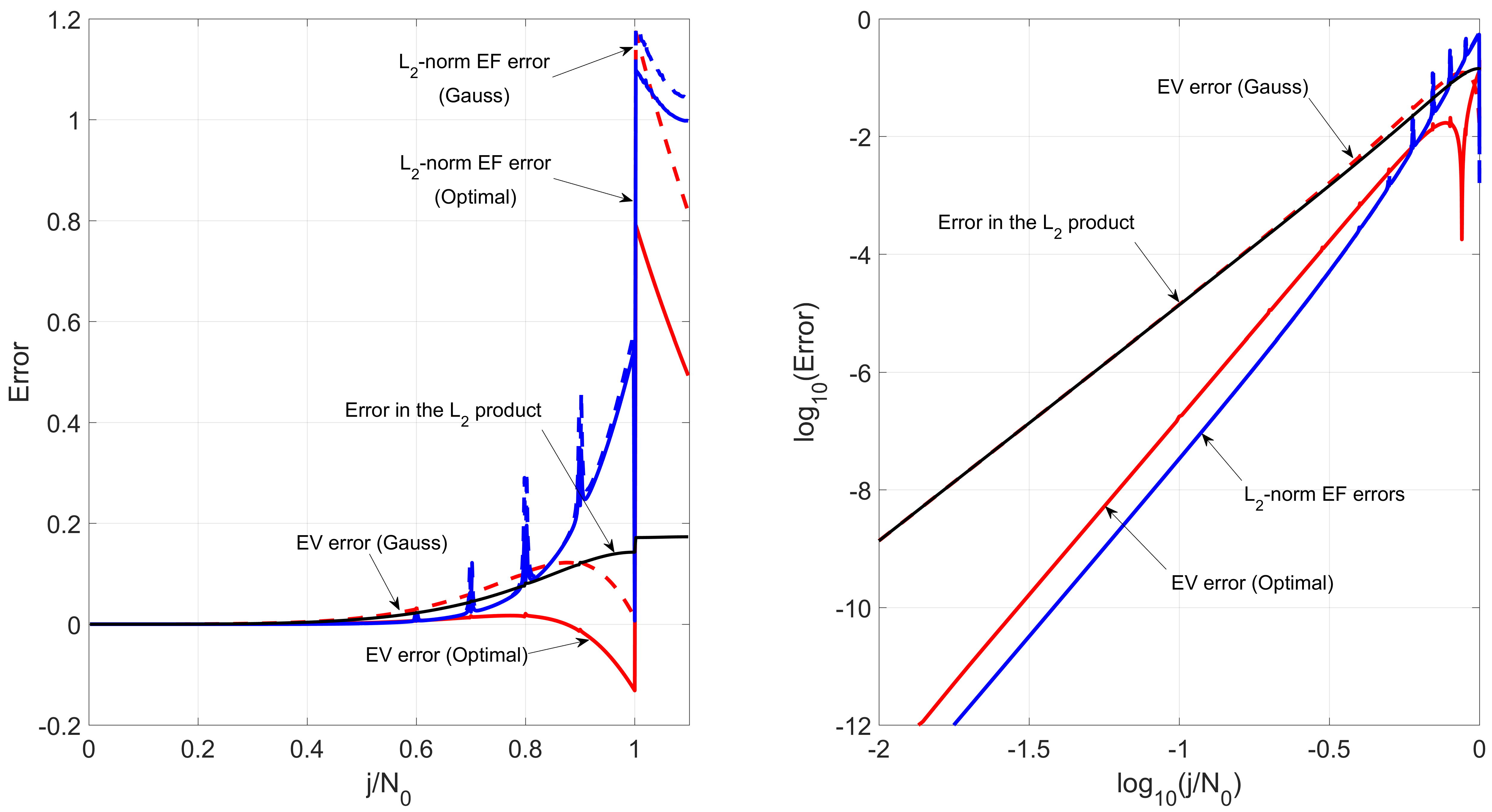}
\caption{Comparison of the eigenvalue and eigenfunction errors of quadratic rIGA using the standard Gauss quadrature (dashed lines) versus the optimal quadrature rule (solid lines). The solid black line is the error in the $L_2$ norm $1 - \left\| {v_j^h} \right\|^2$, that is the last term of \eqref{eq:ptfinal}.}
\end{figure}

Figure 11 compares the accuracy of the quadratic rIGA with the standard Gauss and the optimal quadrature rules. The use of the optimal quadrature results in several orders of magnitude improvement in the eigenvalue errors and increases the convergence rate to $2p+2$. The eigenfunction errors slightly improve in the high-frequency range of the spectrum. The optimal quadratures provide the best approximation properties for the modes in the well-resolved part of the spectrum (many points per wavelength, $j/N_0 \rightarrow 0$). Alternative quadratures can be used for better approximation of the eigenvalues in the middle part of the spectrum, i.e., for practical values of the wavenumber (frequency) \citep{puzyrev2017dispersion}.

\begin{figure}[!ht]
\centering\includegraphics[width=1.0\linewidth]{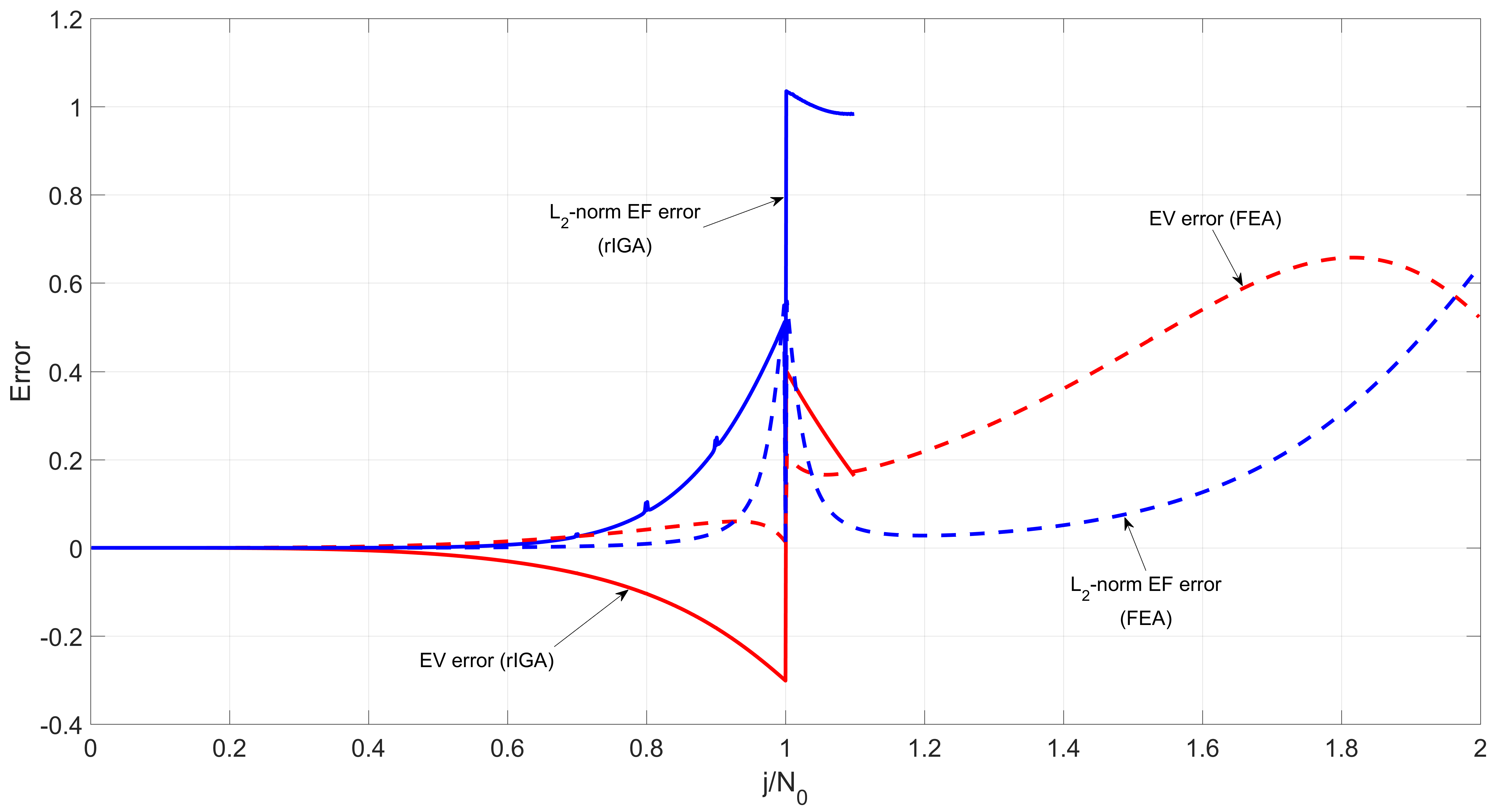}
\caption{Eigenvalue and eigenfunction errors of quadratic rIGA using a non-convex blending of -0.8 Gauss and 1.8 Lobatto quadrature rules (solid lines) versus the standard quadratic $C^0$ elements with the Gauss quadrature (dashed lines).}
\end{figure}

Figure 12 shows the spectra of the quadratic rIGA using a non-convex blending of Gauss and Lobatto quadratures (blending coefficients -0.8 and 1.8, respectively). Using this blending, we decrease the eigenvalue error of the outliers while not degrading the errors in the rest of the spectrum. For comparison purposes, we also show in Figure 12 the errors of the standard, fully integrated FEA for roughly the same number of elements. The eigenvalue errors at the outliers of the quadratic rIGA with this non-convex blending are smaller compared to the optical branch of the FEA thus leading to better time-stepping stability. The convergence rate is the same for both methods, that is $2p$.

Figure 13 compares the errors of the 2D eigenvalue problem using the standard quadratic IGA, rIGA, and FEA versus their optimal counterparts. The IGA and rIGA errors are very similar in the main part of the spectrum; rIGA errors exceed 50\% at the outliers. The quality of the standard FEA approximation is considerably worse. The use of the optimal quadrature leads to large improvement in the accuracy for all methods. Similar to the 1D case shown in Figure 11, optimal rIGA has smaller errors in the high-frequency range of the spectrum. Special quadrature rules may exist that further reduce the impact of the outliers, while delivering an improved convergence rate. For example, a quadrature that significantly underintegrates the basis functions associated with the elements near the boundary and separators may reduce the outlier errors, while not affecting the approximations elsewhere. This topic will be the subject of the future studies.

\begin{figure}[!ht]
\centering\includegraphics[width=1.0\linewidth]{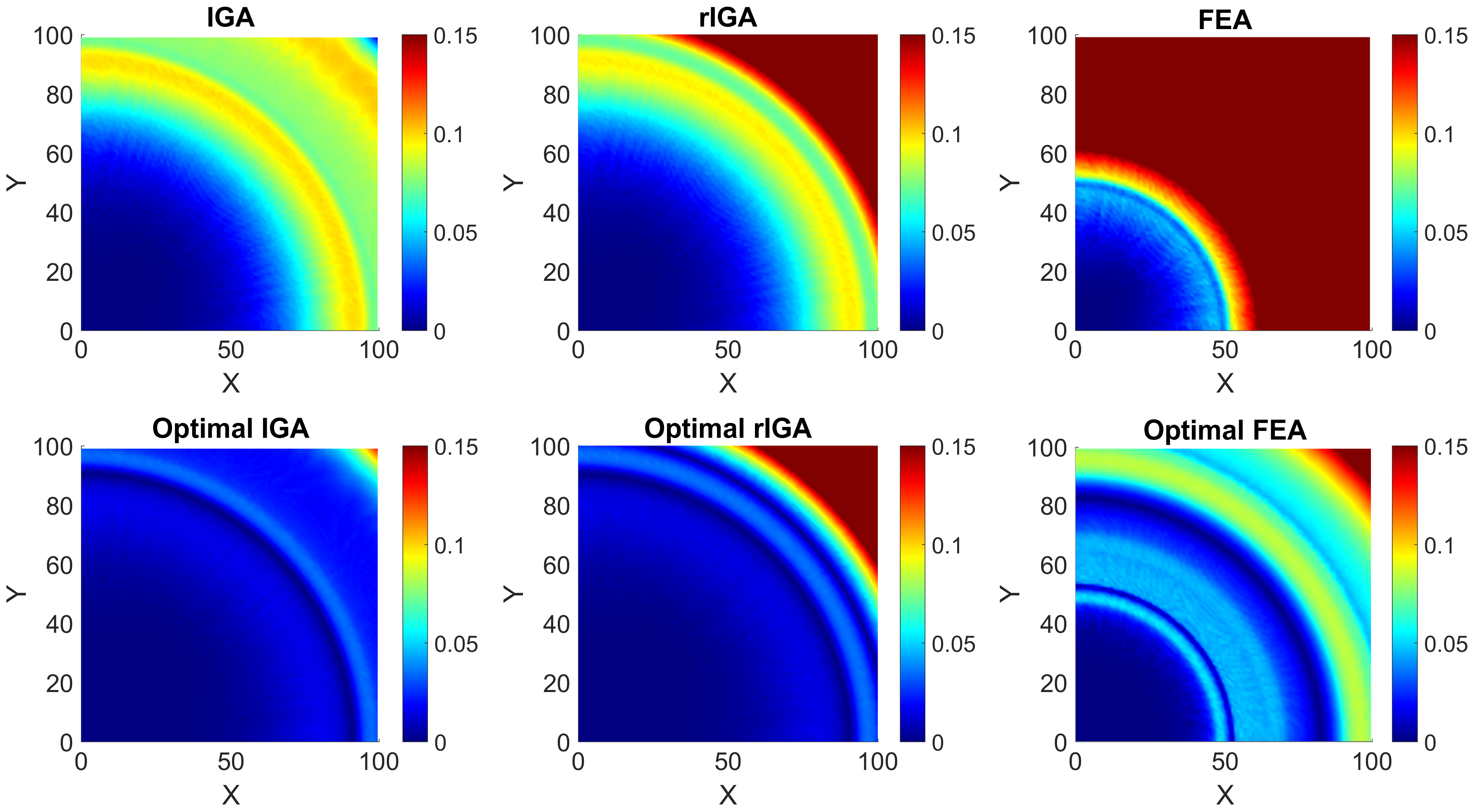}
\caption{Standard quadratic IGA, rIGA with blocks of ten elements, and FEA (top row) versus their optimal counterparts (bottom row). The color shows the relative eigenvalue error.}
\end{figure}





\section{Conclusions and future outlook}

The control of the basis functions continuity is a powerful tool in the \textit{hpk}-refinement space of isogeometric analysis. Highly-continuous isogeometric elements are significantly more accurate than the classical finite elements for the same number of degrees of freedom, though at the cost of an increased computational complexity. One way to overcome this difficulty, is to use isogeometric elements with variable continuity. The resulting refined isogeometric analysis uses blocks of maximum continuity separated by hyperplanes of lower continuity in the mesh to facilitate the elimination of the degrees of freedom by direct solvers. This solution methodology is more than an order of magnitude faster than maximum continuity isogeometric analysis in 3D and several orders of magnitude faster than classical finite elements on a fixed mesh.

We study the spectral approximation properties of refined isogeometric analysis and show how the breaks in continuity and inhomogeneity of the basis affect the errors in the eigenvalues and eigenfunctions. Using properly designed alternative quadratures leads to more accurate results in a similar way as it does in classical isogeometric analysis and finite element methods. Optimal quadrature rules largely reduce the phase error of the method without affecting its overall approximation quality or solution efficiency.

\section{Acknowledgments}

This publication was made possible in part by the CSIRO Professorial Chair in Computational Geoscience at Curtin University and the Deep Earth Imaging Enterprise Future Science Platforms of the Commonwealth Scientific Industrial Research Organisation, CSIRO, of Australia. Additional support was provided by the European Union's Horizon 2020 Research and Innovation Program of the Marie Sklodowska-Curie grant agreement No. 644202. The J. Tinsley Oden Faculty Fellowship Research Program at the Institute for Computational Engineering and Sciences (ICES) of the University of Texas at Austin has partially supported the visits of VMC to ICES. The Spring 2016 Trimester on ``Numerical methods for PDEs'', organised with the collaboration of the Centre Emile Borel at the Institut Henri Poincare in Paris supported VMC's visit to IHP in October 2016.


\nocite{*}

\bibliographystyle{elsarticle-harv}\biboptions{square,sort,comma,numbers}
\bibliography{references}








\end{document}